\documentclass[a4paper]{amsart}
\usepackage[all]{xy}
\usepackage[latin1]{inputenc}
\usepackage{enumerate}
\usepackage[colorlinks=true]{hyperref}
\usepackage{amssymb}

\newtheorem{theorem}{Theorem}[section]
\newtheorem{proposition}[theorem]{Proposition}
\newtheorem{lemma}[theorem]{Lemma}
\newtheorem{corollary}[theorem]{Corollary}
\newtheorem{problem}[theorem]{Problem}

\theoremstyle{definition}
\newtheorem{definition}[theorem]{Definition}
\newtheorem{example}[theorem]{Example}

\theoremstyle{remark}
\newtheorem{remark}[theorem]{Remark}

\numberwithin{equation}{section}

\newcommand{\comment}[1]{}

\newcommand{\Rr}{\mathbb R}

\newcommand{\Nn}{\mathbb N}

\newcommand{\set}[1]{\left\{#1\right\}}

\newcommand{\X}{\ensuremath{\mathfrak{X}}}
\newcommand{\F}{\ensuremath{\mathcal{F}}}
\renewcommand{\d}{\mathrm d}

\DeclareMathOperator{\Der}{Der}

\DeclareMathOperator{\gl}{\mathfrak{gl}}
\DeclareMathOperator{\GL}{GL}

\newcommand{\Jet}{\text{\rm J}}

\newcommand{\G}{\mathcal{G}}            
\newcommand{\s}{\mathbf{s}}             
\renewcommand{\t}{\mathbf{t}}           

\newcommand{\A}{A}                      
\newcommand{\al}{\alpha}                
\newcommand{\be}{\beta}                 
\newcommand{\ga}{\gamma}                
\renewcommand{\gg}{\mathfrak{g}}        
\newcommand{\tto}{\rightrightarrows}    
\newcommand{\wmc}{\omega_{\textrm{MC}}}   

\newcommand{\B}{\mathcal{B}}                    
\newcommand{\BG}{\mathcal{B}_G}                 

\begin{document}

\title{Lie Algebroids and Classification Problems in Geometry}

\author{Rui Loja Fernandes}
\address{Departamento de Matem\'{a}tica,
Instituto Superior T\'ecnico, 1049-001 Lisboa, PORTUGAL}
\email{rfern@math.ist.utl.pt}

\author{Ivan Struchiner}
\address{Instituto de Matem\'{a}tica,
Estat\'{\i}stica e Computa\c{c}\~ao Cient\'{i}fica, Universidade Estadual de
Campinas, BRASIL} \email{ivan@ime.unicamp.br}

\thanks{The first author was supported in part by FCT/POCTI/FEDER and by grants POCI/MAT/55958/2004 and POCI/MAT/57888/2004.
The second author was supported in part by FAPESP grant 03/13114-2
and CAPES grant BEX3035/05-0.}

\begin{abstract}
We show how one can associate to a given class of finite type $G$-structures a classifying Lie algebroid.
The corresponding Lie groupoid gives models for the different geometries that one can find in the class,
and encodes also the different types of symmetry groups.
\end{abstract}

\maketitle

\tableofcontents

\section{Introduction}

The main purpose of this paper is to describe how one can associated
to certain classes of geometric structures a \emph{classifying Lie algebroid} $A \rightarrow X$,
which has the following properties:
\begin{enumerate}[(i)]
\item to each point on the base $X$ there corresponds a (germ of a) geometric
structure of the class;
\item two structures in the class are locally isomorphic if and only if they correspond
to the same point of $X$;
\item the isotropy Lie algebra at a point is the symmetry Lie
algebra of the corresponding geometric structure;
\item two points belong to the same orbit of $A$ if and only if there exists a geometric structure of the class which contains their corresponding germs.
\end{enumerate}
We will also be interested in the associated Lie groupoid which allows one
to construct different models for the distinct geometries that one can find in the class.

The geometric structures that will be considered here are $G$-structures. The notion of
$G$-structure is quite general and includes most of the classical geometric structures.
Since we want to avoid technical difficulties with infinite dimensional spaces, we restrict
ourselves to $G$-structures of finite type, which means that $G$ is a Lie group of
finite type (see Definition \ref{finitetype}). For finite type $G$-structures the
classification problem can be reduced to a classification problem of $\set{e}$-structures
(i.e., coframes) through the method of prolongation.

The method then goes as follows. Suppose that through prolongation
we have reduced our classification problem to a problem for a
class of coframes. Under suitable regularity assumptions, any
coframe is determined by a finite set of functionally independent
\emph{structure invariants} (this is the generic case). It will
then follow that a regular coframe determines a \emph{transitive}
Lie algebroid. By a \emph{class of coframes}, we mean a family of
coframes in which all members are determined by the same (not
necessarily independent) structure invariants. In order to obtain
a classifying Lie algebroid over a finite dimensional base $X$, we
must restrict ourselves to classes of coframes determined by a
finite set of structure invariants.

Finally, we will be interested in finding explicit models for the different
geometric structures belonging to the classes being studied. The problem of
finding a coframe with given structure invariants is known as
\emph{Cartan's realization problem} (\cite{Bryant,Cartan}).
We will prove the following result:

\begin{theorem}
Given the initial data of a Cartan's realization problem, one has that:
\begin{enumerate}[(i)]
\item a realization exists if and only if the initial data form the structure
functions of a Lie algebroid $A$;
\item any realization is locally equivalent to a neighborhood
of the identity of an $\mathbf{s}$-fiber of a (local) Lie groupoid
$\mathcal{G}$ integrating $A$, equipped with the Maurer-Cartan form;
\item any two germs of coframes belong to the same global coframe if and only if they
correspond to points on the same orbit of $A$.
\end{enumerate}
\end{theorem}

We call the Lie algebroid in the theorem the \emph{classifying Lie algebroid}
of the class of coframes. The theorem then shows that, when this Lie algebroid
is integrable (see \cite{CrainicFernandes}), each $\s$-fiber of its groupoid,
equipped with its Maurer-Cartan form furnishes an explicit \emph{universal}
model of the coframe.

This paper is a preliminary announcement of the results obtained by the second author
in his thesis. In particular, a detailed analysis of specific
examples is beyond the scope of this paper.

This paper is organized as follows. In Section 2, we recall some classical results
concerning finite type $G$-structures that we will need, including the method of
prolongation of $G$-structures. In Section 3, which contains the main results, we
start by constructing the classifying Lie algebroid for $\{e\}$-structures. Then we
explain Cartan's realization problem for these structures, and answer the two basic
questions: (i) existence and (ii) classification of realizations. Along the way, we
give a brief study of Maurer-Cartan forms on Lie groupoids, since these play a crucial
role in the classification. We end this section explaining how to extend these results
for any finite type $G$-structures. In Section 4, we illustrate briefly our approach
with two examples: constant curvature Riemannian metrics and Bochner-K\"{a}hler metrics.
\vskip 10 pt

\noindent \textbf{Acknowledgments.} The first author would like to
thank Professor Luiz San Martin for its hospitality during 
a visit to UNICAMP. The second author would like to thank the Center for Mathematical
Analysis, Geometry and Dynamical Systems of IST, for its support
during a visit where this work was initiated, and also Professor Luiz San Martin 
for all the helpful discussions and advice.

\section{Finite Type $G$-Structures}

In this section we recall the basic facts from the theory of $G$-structures that we will use.
We refer to \cite{Kobayashi,SingerSternberg,Sternberg} for details.

\subsection{$G$-Structures}

Denote by
\[\xymatrix{\B(M) \ar[d]_{\pi} & \ar@(dr,ur)@<-5ex>[]_{\GL_n}\\
M}
\]
the bundle of frames on $M$. This principal bundle carries a canonical
1-form with values in $\Rr^n$, denoted by $\omega\in\Omega^{1}(\B(M);\Rr^{n})$,
and which is defined by
\[ \omega_{p}(X):= p^{-1}(\pi_{\ast}X),\quad (X\in T_{p}\mathcal{B}(M)).\]
This form is called the \textbf{tautological form} (or \emph{soldering form}) of
$\mathcal{B}(M)$. It is a tensorial form, i.e, it is horizontal and $\GL(n)$-equivariant
(with respect to the defining action on $\Rr^n$). Note that a subspace $H_p\subset T_p\BG(M)$
is horizontal iff the restriction $\omega:H_p\to\Rr^n$ is an isomorphism.

Every diffeomorphism $\varphi$ between two manifolds $M$ and $N$ lifts to an isomorphism of
the associated frame bundles:
\[\B(\varphi): \B(M) \rightarrow \B(N).
\]
The correspondence which associates to each manifold its frame
bundle and to each diffeomorphism its lift is functorial.

Now let $G$ be a Lie subgroup of $\GL(n)$. Recall that a \textbf{$G$-structure}
is a reduction of the frame bundle $\B(M)$ to a principal $G$-bundle. This means that
$\B_G(M) \subset \B(M)$ is a sub-bundle such that for any $p \in \B_G(M)$
and $a \in \GL(n)$ we have $pa \in \B_G(M)$ if and only if $a \in G$. Given a $G$-structure
$\B_{G}(M)$, we will still denote by $\omega$ the restriction of the tautological form to
$\B_{G}(M)$.

\begin{definition}
Two $G$-structures $\BG(M)$ and $\BG(N)$ are said to be
\textbf{equivalent} if there exists an diffeomorphism $\varphi:M
\rightarrow N$ such that
\[\B(\varphi)(\BG(M)) = \BG(N)
\]
\end{definition}

\subsection{Equivalence of $G$-Structures}

One of the basic problems we will be interested is deciding if two $G$-structures
are equivalent. The tautological form is the clue to the solution of this
equivalence problem. The reason is the following result:

\begin{proposition}
Two $G$-structures over $M$ and $N$ are equivalent if and only if
there exists a principal $G$-bundle isomorphism $\psi : \B_{G}(M) \rightarrow
\B_{G}( N)$ such that $\psi^{\ast}\omega_{N} = \omega_{M}$.
\end{proposition}

In order to obtain an invariant of equivalence of $G$-structures, let us
choose some horizontal space $H_p$ at $p\in \BG(M)$. Given $v\in\Rr^n$ there exists a unique
$\widetilde{v}\in H_p$ such that $\omega(\widetilde{v})=v$, so one defines:
\begin{align}\label{firetorderstructurefunction}
c_{H_p}:\wedge^{2} \Rr^{n} & \to \Rr^{n}, \\
c_{H_p}(v,w)&:=\d\omega(\widetilde{v},\widetilde{w}).\nonumber
\end{align}
This depends on the choice of horizontal space, so it does not define an invariant yet.
If $H_p$ and $H'_p$ are two distinct horizontal spaces at $p\in \BG(M)$, then one checks that:
\[ c_{H_p}-c_{H'_p}\in \mathcal{A}(\hom(\Rr^n,\gg)),\]
where $\gg\subset\gl(n)$ is the Lie algebra of $G$ and $\mathcal{A}$ denotes the
anti-symmetrization operator:
\begin{align*}
\mathcal{A}:\hom(\Rr^{n},\mathfrak{g}) & \rightarrow \hom(\wedge^{2}\Rr^{n},\Rr^{n}), \\
\mathcal{A}(T)(u,v) & := T(u)v - T(v)u.
\end{align*}
Hence, we can set:

\begin{definition}
Given a $G$-structure $\B_G(M)$ one defines its \textbf{first order structure function}:
\[
c:\B_G(M)\rightarrow\frac{\hom(\wedge^{2}\Rr^n,\Rr^n)}{\mathcal{A}(\hom(\Rr^n,\gg))},\quad
c(p):=[c_{H_p}].
\]
\end{definition}

Since an isomorphism $\psi: \B_{G}(M)\rightarrow \B_{G}(N)$ maps horizontal spaces to
horizontal spaces and it is an equivalence if and only if $\psi^{\ast}\omega_{N} = \omega_{M}$,
we see that

\begin{proposition}
Let $\B_{G}(M)$ and $\B_{G}(N)$ be $G$-structures. If $\phi : M
\rightarrow N$ is an equivalence then
\[
c_{N} \circ \mathcal{B}(\phi) = c_{M}.
\]
\end{proposition}

\subsection{Prolongation}
In order to obtain more refined invariants of equivalence of $G$-structures
one needs to look at higher order terms. This process is known as \emph{prolongation}
and takes place on the jet bundles $\Jet^k\BG(M)$.

Let $\pi:E\to M$ be a fiber bundle. We denote by $\pi^1:\Jet^1E\to M$ its first
jet bundle, which has fiber over $x\in M$:
\[ (\Jet^1 E)_x=\set{j^1_x s|\text{ $s$ a section of $E$}}.\]
This bundle can also be described geometrically as:
\[ \Jet^1E=\set{H_p: p\in E\text{ and }H_p\subset T_p E\text{ horizontal}}.\]
If one defines the projection $\pi^1_0:\Jet^1E\to E$ by $\pi^1_0(H_p)=p$, then
$\Jet^1E$ is an affine bundle over $E$.

\begin{example}
For a $G$-structure $\BG(M)$ the first structure function can also be described
as a function $c:\Jet^1\BG(M)\to \hom(\wedge^{2}\Rr^n,\Rr^n)$ by formula (\ref{firetorderstructurefunction}).
\end{example}

It is easy to see that in the case of the frame bundle
$\pi:\B(M)\to M$ its first jet bundle
$\pi^1_0:\Jet^1\B(M)\to\B(M)$ can be identified with a sub bundle
of $\B(\B(M))$: to a horizontal space $H_p\subset T_p\B(M)$ we
associate a frame in $\B(M)$ (which we view as an isomorphism
$\phi:\Rr^n\times \gl(n)\to T_p\B(M)$):
\[ \Rr^n\times \gl(n)\ni(v,\xi)\stackrel{\phi}{\longmapsto} (\pi|_{H_p}\circ p)^{-1}(v)+\xi\cdot p\in T_p\B(M).\]
Note that if $H_p$ and $H'_p$ are two horizontal spaces at $p\in\B(M)$, the corresponding frames
$\phi,\phi':\Rr^n\times \gl(n)\to T_p\B(M)$ are related by:
\[ \phi'(v,\xi)=\phi(v,\xi)+T(v)\cdot p,\]
for some $T\in\hom(\Rr^n,\gl(n))$. Conversely, given a frame $\phi$ associated with some horizontal space $H_p$ and
$T\in\hom(\Rr^n,\gl(n))$, this formula determines a frame $\phi'$ which is associated with another
horizontal space $H'_p$. It follows that $\Jet^1\B(M)$ is a $\hom(\Rr^n,\gl(n))$-structure,
where we view $\hom(\Rr^n,\gl(n))\subset\GL(\Rr^n\oplus\gl(n))$ as the subgroup formed by those
transformations:
\[ (v,\xi)\mapsto (v,\xi+T(v)),\text{ with }T\in\hom(\Rr^n,\gl(n)).\]

Assume now that $\BG(M)$ is a $G$-structure so that $J^1\BG(M)\subset J^1\B(M)$ is a sub-bundle. An argument
entirely similar to one just sketched gives:

\begin{proposition}
If $\BG(M)$ is a $G$-structure then $J^1\BG(M)\to\BG(M)$ is a $\hom(\Rr^n,\gg)$-structure.
\end{proposition}

In order to motivate our next definition we look at a simple example.

\begin{example}
Let us consider the flat $G$-structure on $\Rr^n$:
\[ \BG(\Rr^n):=\Rr^n\times G\subset \B(\Rr^n)=\Rr^n\times \GL(n).\]
Given a vector field $X$ denote by $\phi^t_X:\Rr^n\to \Rr^n$ its flow. Observe
that $X$ is an infinitesimal automorphism of the $G$-structure $\BG(\Rr^n)$
iff $\phi^t_X$  lifts to an automorphism $\mathcal{B}(\phi^t_X):\BG(\Rr^n)\to\BG(\Rr^n)$.
The lifted flow $\mathcal{B}(\phi^t_X)$ is the flow of a lifted vector field on $\BG(\Rr^n)$:
in coordinates $(x^1,\dots,x^n)$, so that $X=X^i\frac{\partial}{\partial x^i}$,
the lifted vector field is given by:(\footnote{We use the convention of summing over repeated indices.})
\[ \widetilde{X}=\frac{\partial X^i}{\partial x^j}\frac{\partial }{\partial p^i_j},\]
where $(p^i_j)$ are the associated coordinates in $\B(\Rr^n)$ so that a frame $p\in\B(\Rr^n)$
is written as:
\[ p=(p^i_1\frac{\partial }{\partial x^i},\dots,p^i_n\frac{\partial }{\partial x^i}).\]
It follows that $X$ is an infinitesimal automorphism iff:
\[ \left[\frac{\partial X^i}{\partial x^j} \right]_{i,j=1,\dots,n}\in\gg\subset\gl(n).\]
Let us assume now that the lifted flow fixes $(0,I)\in\Rr^n\times\GL(n)$. The lifted vector field
$\widetilde{X}$ vanishes at this point. If we now prolong to the jet bundle $J^1\BG(\Rr^n)$,
we obtain a flow which is generated by a vector field:
\[ j^1\widetilde{X}=\frac{\partial X^i}{\partial x^{j_1}\partial x^{j_2}}\frac{\partial }{\partial p^i_{j_1,j_2}},\]
where $(x^i,p^i_j,p^i_{j_1,j_2})$ are the induced coordinates on the jet bundle. Note that the coefficients
$a^i_{j_1,j_2}=\frac{\partial X^i}{\partial x^{j_1}\partial x^{j_2}}$ of $j^1\widetilde{X}$ satisfy:
\[ \left[a^i_{j_1,j_2}\right]_{i,j_1=1,\dots,n}\in\gg\subset\gl(n),\]
and are symmetric in the indices $j_1$ and $j_2$. Hence, we conclude that:

\begin{lemma}
The lifts of the symmetries of the flat $G$-structure $\BG(\Rr^n)=\Rr^n\times G$ to the
jet space $\Jet^1\BG(\Rr^n)$ generate a Lie subgroup $G^{(1)}\subset\hom(\Rr^n,\gg)$
with Lie algebra:
\[ \gg^{(1)}:=\set{T\in\hom(\Rr^n,\gg):T(u)v=T(v)u,\forall u,v\in\Rr^n}.\]
\end{lemma}
\end{example}

This motivates the following definition:

\begin{definition}\label{finitetype}
Let $\gg\subset\gl(V)$ be a Lie algebra. The \textbf{first prolongation} of
$\gg$ is the subspace $\gg^{(1)}\subset\hom(V,\gg)$
consisting of those $T:V\to\gg$ such that
\[
T(v_1)v_2 = T(v_1)v_2,\quad \forall v_1,v_2\in V.
\]
The \textbf{$k$-th prolongation} of $\gg$ is the subspace $\gg\subset\hom(V,\gg^{(k-1)})$
defined inductively by
\[
\gg^{(k)} = (\gg^{(k-1)})^{(1)}.
\]
A Lie algebra $\gg$ is said to be of finite type $k$ if there exists $k \in
\mathbb{N}$ such that $\gg^{(k-1)}\ne 0$ and $\gg^{(k)}= 0$.
\end{definition}

Similarly, at the group level, one introduces:

\begin{definition}
Let $G$ be a subgroup of $GL(V)$. The \textbf{first prolongation} of $G$
is the subgroup $G^{(1)}$ of $GL(V\oplus \gg)$ consisting of those
transformations of the form:
\[ (v,\xi)\mapsto (v,\xi+T(v)),\text{ with }T\in\gg^{(1)}.\]
Similarly, the \textbf{$k$-th prolongation} of $G$ is the subgroup $G^{(k)}$
of $\GL(V\oplus\gg\oplus\gg^{(1)}\oplus\cdots\oplus\gg^{(k)})$ defined inductively by
\[
G^{(k)}:=(G^{(k-1)})^{(1)}.
\]
\end{definition}

Note that the prolongations $G^{(k)}$ are all Abelian groups. Now, to each $G$-structure
$\BG(M)$ we can always reduce the structure group of $\Jet^1\BG(M)$ to $G^{(1)}$,
obtaining a $G^{(1)}$-structure:

\begin{proposition}
Let $\BG(M)$ be a $G$-structure over $M$ with first structure function $c:\Jet^1\BG(M)\to \hom(\wedge^{2}\Rr^n,\Rr^n)$.
Each choice of a complement $C$ to $\mathcal{A}(\hom(\Rr^n,\gg))$ in $\hom(\wedge^{2}\Rr^n,\Rr^n)$ determines
a sub-bundle:
\[ \BG(M)^{(1)}=\set{H_p\in \Jet^1\BG(M):c_{H_p}\in C},\]
which is a reduction of $\Jet^1\BG(M)$ with structure group
$G^{(1)}$. Different choices of complements determine sub-bundles
which are related through right translation by an element in
$\hom(\Rr^n,\gg)$.
\end{proposition}

The $G^{(1)}$-structure $\BG(M)^{(1)}\to \BG(M)$ is called the
\textbf{first prolongation} of $\BG(M)$. Similarly, working
inductively, one defines the \textbf{$k$-th prolongation} of
$\BG(M)$:
\[ \BG(M)^{(k)}=(\BG(M)^{(k-1)})^{(1)},\]
which is $G^{(k)}$-structure over $\BG(M)^{(k-1)}$.

The relevance of prolongation for the problem of equivalence is justified
by the following basic result:

\begin{theorem}
Let $\BG(M)$ and $\BG(N)$ be $G$-structures. Then
$\BG(M)$ and $\BG(N)$ are equivalent if and only if their
first prolongations $\BG(M)^{(1)}$ and $\BG(N)^{(1)}$ are equivalent
$G^{(1)}$-structures.
\end{theorem}

One can now obtain new necessary conditions for equivalence by
looking at the structure function of the prolongation
$\BG(M)^{(1)}$ which is a a function
\[
c^{(1)} : \BG(M)^{(1)} \rightarrow \frac{\hom(\wedge^{2}(\Rr^{n} \oplus \gg),\Rr^{n} \oplus
\gg)}{\mathcal{A}(\hom(\Rr^{n} \oplus \gg, \gg^{(1)}))}
\]
called the \textbf{second order structure function} of $\BG(M)$.
Then one can continue this process by constructing the second
prolongation and analyzing it's structure function and so on.

Thus, the importance of structures of finite type is that we can
reduce the set of necessary conditions for checking that two
$G$-structures are equivalent to a finite amount. In fact, by the
method of prolongation, the equivalence problem for finite type
$G$-structures reduces to an equivalence problem for $\{e\}$-structures
(coframes). Moreover, one can show that $G$-structures of finite type
always have finite dimensional symmetry groups.

\subsection{Second Order Structure Functions}
\label{sec:2ndorder}

By working inductively, all we really must understand are the
second order structure functions which we now describe.

Let $z =H_p \in \BG^{(1)}=\BG(M)^{(1)}$ and let
$\mathcal{H}_{z}$ be a horizontal subspace of $T_z\BG^{(1)}$.
Then $c_{\mathcal{H}_z}^{(1)}\in\hom(\wedge^{2}(\Rr^n
\oplus \gg), \Rr^n \oplus \gg)$ and we decompose it into three
components:
\begin{align*}
\hom(\wedge^{2}(\Rr^n\oplus \gg),\Rr^n \oplus \gg) =
\hom(\wedge^{2}\Rr,&\Rr^n\oplus \gg)\oplus\hom(\Rr^n\otimes\gg,\Rr^n\oplus \gg)\oplus\\
&\oplus\hom(\wedge^{2}\gg,\Rr^n\oplus \gg)
\end{align*}
Let us describe each of the components of the (representative of the)
second order structure function. We denote by $u,v$ elements
of $\Rr^n$ and by $A,B$ elements of $\gg$:
\begin{itemize}
\item The first component of $c_{\mathcal{H}_{z}}^{(1)}$ includes the structure
function of $\BG(M)$:
\[
c_{\mathcal{H}_{z}}^{(1)}(u,v) = c_{H_{p}}(u,v) +
b_{\mathcal{H}_{z}}(u,v),
\]
for some $b_{\mathcal{H}_{z}} \in \hom(\wedge^{2}\Rr^{n},\gg)$.

\item The second component of $c_{\mathcal{H}_{z}}^{(1)}$ has the form:
\[
c_{\mathcal{H}_{z}}^{(1)}(A,u) = -Au + S_{\mathcal{H}_{z}}(A,u)
\]
for some $S_{\mathcal{H}_{z}} \in \hom(\Rr^n\otimes\gg,\Rr^n)$.

\item The last component of $c_{\mathcal{H}_{z}}^{(1)}$ is
given by
\[
c_{\mathcal{H}_{z}}^{(1)}(A,B) = -[A,B]_{\gg} \text{.}
\]
\end{itemize}

An important special case occurs when $G^{(1)} = \{e\}$. In this case,
a $G^{(1)}$-structure amounts to choosing a horizontal space at each
$p\in\BG(M)$, which in turn is the same as picking a $\gg$-valued (not necessarily
equivariant) form $\phi$ on $\BG$. The pair $(\omega,\phi)$ is a coframe on $\B_G$.
Now, in this case, the projection from $\B_{G^{(1)}}$ onto $\B_G$ is a
diffeomorphism, so we may view the second order structure functions
as functions on $\B_G$. If we do this, we obtain the \textbf{structure
equations} of the pair $(\omega,\phi)$:
\begin{equation}
\label{eq:structure}
\left\{
\begin{array}{l}
\d \omega=c \circ \omega \wedge \omega - \phi \wedge \omega\\
\\
\d \phi= b \circ \omega \wedge \omega + S \circ \phi \wedge
\omega - \phi\wedge\phi
\end{array}
\right.
\end{equation}
where $\phi\wedge\omega$ is the $\Rr^n$-valued 2-form obtained from the
$\gg$-action on $\Rr^n$ and $\phi \wedge \phi$ is the
$\gg$-valued 2-form obtained from the Lie bracket on $\gg$.

If, additionally, the horizontal spaces can be chosen right invariant (so that ${R_a}_{\ast}H_{p}=H_{pa}$ for
all $a\in G$), we obtain a principal bundle connection on $\BG$ with connection form
$\phi$. In this case, we find that:
\begin{itemize}
\item $S$ vanishes identically;
\item $b$ is the curvature of the connection;
\end{itemize}
so we see that, in this case, equations (\ref{eq:structure}) reduce to the usual structure equations for a connection.

\section{Cartan's Realization Problem}

Given a reasonable class of $G$-structures we now explain how one
can associate to it a \emph{classifying Lie algebroid}.

\subsection{Equivalence of Coframes}
Assume that we have prolonged our finite type $G$-structure as
much as necessary so we arrive at an $\set{e}$-structure. Since an
$\set{e}$-structure is just the specification of a frame (or a
coframe), we must then solve a problem of equivalence of coframes.
We recall here how this can be dealt with. For details we refer to
\cite{Olver} and \cite{Sternberg}.

Let $\theta=\{\theta^{1},...,\theta^{n}\}$ be a coframe (i.e., a spanning set of
everywhere linearly independent $1$-forms) on an $n$-dimensional manifold $M$.
If $\bar{M}$ is another $n$-dimensional manifold and $\bar{\theta}=\{\bar{\theta}^{i}\}$ a coframe on
$\bar{M}$, the equivalence problem asks:(\footnote{We will use unbarred letters to
denote objects on $M$ and barred letters to denote objects on $\bar{M}$.})

\begin{problem}[Equivalence Problem]
\label{PEquivalence} Does there exist a (locally defined)
diffeomorphism $\phi:M\rightarrow\bar{M}$ satisfying
\[
\phi^{\ast}\bar{\theta}^{i}=\theta^{i}\text{?}
\]
\end{problem}

Exterior differentiation of the 1-forms $\theta^k$ give the \textbf{structure equations}:
\begin{equation}
\label{structureeqcoframe}
\d\theta^{k}=\sum_{i<j}C_{ij}^{k}(x)\theta^{i}\wedge\theta^{j},
\end{equation}
for some functions $C_{ij}^{k}\in\mathrm{C}^{\infty}(M)$ called the \textbf{structure functions}
of the coframe. These functions play a crucial role in the
study of the equivalence problem. For example, since for any coframe
$\bar{\theta}$ equivalent to $\theta$ one must have
\[
\bar{C}_{ij}^{k}(\phi(x))=C_{ij}^{k}(x),
\]
the structure functions furnish a a set of invariants
of the equivalence problem.

\begin{definition}
A function $I\in\mathrm{C}^{\infty}(M)$ is called an invariant
function of a coframe $\{\theta^{i}\}$ if for any locally defined
self equivalence (symmetry) $\phi:M\to M$ one has
\[
I\circ\phi=I.
\]
\end{definition}

Now, for any function $f\in\mathrm{C}^{\infty}(M)$ one defines
its coframe derivatives $\frac{\partial f}{\partial\theta^{k}}$
as the coefficients of the differential of $f$ when
expressed in terms of the coframe $\{\theta^{i}\}$,
\[
\d f=\sum_{k}\frac{\partial f}{\partial\theta^{k}}\theta^{k}.
\]
Using the fact that $\d \phi^{\ast} = \phi^{\ast} \d$, it
follows that if $I\in\mathrm{C}^{\infty}(M)$ is an invariant
function, then so are its coframe derivatives $\frac{\partial I}{\partial\theta^{k}}$,
for all $1\leq k\leq n$. It is then natural to consider the sets of
\emph{structure invariants},
\[
\mathcal{F}_{s} =  \left\{  C_{ij}^{k},\frac{\partial C_{ij}^{k}}
{\partial\theta^{l}},\ldots,\frac{\partial^{s}C_{ij}^{k}}{\partial
\theta^{l_{1}}\cdots\partial\theta^{l_{s}}}\right\}
\]
which give us an infinite number of necessary conditions to solve
the equivalence problem. Schematically, we may write
\begin{equation}
\label{necessary}
\phi^{\ast}\mathcal{\bar{F}}_{s}=\mathcal{F}_{s}
\end{equation}
for all $s \geq 0$.

In order to be able to proceed, we must first reduce these necessary
conditions to a finite number.

\begin{definition}
A coframe $\theta=\{\theta^{i}\}$ is called \textbf{fully regular} if
for each integer $s=0,1,2,\dots$ the $s^\text{th}$ order structure map $C^{(s)}:M\to\Rr^{N_s}$,
whose components are the structure invariants in $\mathcal{F}_{s}$, is
regular (i.e., has constant rank).
\end{definition}

Note that, in the fully regular case, locally we can always find a finite set of
functionally independent structure invariants
$\set{h_1,\ldots,h_d}$ which generate the full set of structure
invariants. This means that for every $s \geq 0$, any element in
$f\in \F_s$ can be written as:
\[
f= H(h_{1},...,h_{d}),
\]
for some function $H:\Rr^d \rightarrow \Rr$. Finally,
observe that in order to verify \eqref{necessary}, it suffices to verify that
\[
\phi^{\ast}{\bar{h}}_{i}= h_{i}
\]
for $1 \leq i \leq d$. The integer $d$ is called the \textbf{rank of the coframe}.

We can now summarize all essential data obtained from a regular
coframe of rank $d$ in the following convenient manner. First
of all, we have a set of invariant functions which determine a map
\[
h:M \rightarrow \Rr^{d},\quad h(x):=(h_{1}(x),...,h_{d}(x)).
\]
Next, since $\{h_{i}\}$ are independent and generate
$\mathcal{F}_{t}$ for all $t \geq 0$ (in particular
$\mathcal{F}_{0}$), we can think of $h_{1},...,h_{d}$ as
coordinates on an open subset $X$ of $\Rr^{d}$. Then the structure
functions may be seen as functions $C_{ij}^{k} \in
\mathrm{C}^{\infty}(X)$. Finally, differentiating $h_{a}$ we
obtain
\[
\d h_{a}=\sum_{i}F_{i}^{a}\theta^{i},
\]
where, for the same reason, $F_{i}^{a}$ are invariant functions
and hence can be seen as elements of $\mathrm{C}^{\infty}(X)$. In
other words, our final structure equations are:
\begin{align*}
\d \theta^{k}&=\sum_{i<j}C_{ij}^{k}(h)\theta^{i}\wedge\theta^{j},\\
\d h_{a} &=\sum_{i}F_{i}^{a}(h)\theta^{i}.
\end{align*}

As we will see next, the functions $F_{i}^{a},C_{ij}^{k}\in
\mathrm{C}^{\infty}(X)$ form the initial data of a \emph{Cartan's
realization problem}. It will then be clear that to any fully
regular coframe we can associate a transitive, flat Lie algebroid.

\subsection{Cartan's Realization Problem}

Not every set of functions $F_{i}^{a},C_{ij}^{k}\in
\mathrm{C}^{\infty}(X)$ determines a class of coframes.
Determining when this is true is the content of:

\begin{problem}[Cartan's Realization Problem]  \label{PRealization}
One is given:
\begin{itemize}
\item an integer $n\in\mathbb{N}$,
\item an open set $X\subset\Rr^{d}$,
\item a set of functions $C_{ij}^{k}\in\mathrm{C}^{\infty}(X)$ with indexes $1\leq i,j,k\leq n$,
\item and a set of functions $F_{i}^{a}\in\mathrm{C}^{\infty}(X)$ with $1\leq a\leq d$
\end{itemize}
and asks for the existence of
\begin{enumerate}[(1)]
\item an $n$-dimensional manifold $M$
\item a coframe $\{\theta^{i}\}$ on $M$
\item a function $h : M \rightarrow X$
\end{enumerate}
satisfying the structure equations
\begin{align}
\d \theta^{k}&=\sum_{i<j}C_{ij}^{k}(h)\theta^{i}\wedge\theta^{j},\label{dtheta}\\
\d h_{a} &=\sum_{i}F_{i}^{a}(h)\theta^{i}.\label{dh}
\end{align}
\end{problem}

A solution $(M,\theta^{i},h)$ to Cartan's problem is called \textbf{realization}. If for
each $h_{0}\in X$ there exists a realization $\left(M,\theta^{i},h\right)$ and
$x_{0}\in M$ such that $h\left(x_{0}\right)=h_{0}$, the initial set of data
is said to specify a \textbf{class of coframes}. The classification
problem can now be stated as:

\begin{problem}[Classification Problem]
\label{PClassification}
What are the possible solutions to a Cartan's realization problem?
When are two solutions to a Cartan's realization problem equivalent?
\end{problem}

\subsection{Existence of Solutions}
Let us consider the question of existence of solutions to a
Cartan's realization problem. One obtains some obvious necessary
conditions for existence by differentiating the structure
equations above and using $\d^{2}=0$. The resulting equations is a
complicated set of non-linear pde's, but which have a very simple
geometric interpretation: they are the differential equations that
define a Lie algebroid. In fact, denoting by $A=X\times \Rr^n\to
X$ the trivial vector bundle over $X$ of rank $n$, with basis of
sections $\{e_1,\dots,e_n\}$, we define a bracket
\begin{equation}
\label{eq:bracket}
[e_i,e_j]_A:=C_{ij}^k e_k,
\end{equation}
and an anchor map $\#:A\to TX$ by:
\begin{equation}
\label{eq:anchor}
\#(e_i)=F_i^a \frac{\partial}{\partial x^a}.
\end{equation}
Then one checks that:
\begin{align*}
\d^2\theta^k=0 &\Longleftrightarrow \text{ Jacobi identity for }[~,~]_A\\
\d^2h_a=0  &\Longleftrightarrow \text{ $\#:\Gamma(A)\to\X(M)$ is a morphism.}
\end{align*}
Thus, we have obtained:

\begin{proposition}
\label{NecessaryRealization} For a Cartan's realization problem to
have solutions for every $h_0 \in X$ it is necessary that
$C_{ij}^{k},F_{i}^{a}\in\mathrm{C}^{\infty}(X)$ be the structure
functions of a flat Lie algebroid $A$ over $X$.
\end{proposition}

It turns out that the above necessary condition is also sufficient, as we shall
show in the sequel.

\subsection{The Maurer-Cartan form on a Lie groupoid}
We will need a generalization of the usual Maurer-Cartan form on a Lie group
to a Lie groupoid.

By a differential $1$-form on a manifold $M$ with values in a Lie algebroid
$A\rightarrow X$ we mean a bundle map
\[
\xymatrix{ TM \ar[d] \ar[r]^{\eta} & A \ar[d] \\
M \ar[r]_h & X }
\]
which is compatible with the anchors, i.e, such that
\[ \xymatrix{
TM \ar[r]^{\eta} \ar[dr]_{h_{\ast}} &  A \ar[d]^{\#} \\
& TX}
\]
Since there is no canonical way of differentiating forms with
values in a vector bundle, we introduce an arbitrary connection
$\nabla$ on $A\rightarrow X$, and for $\xi_{1},\xi_{2} \in \X(M)$
we define
\[ d_{\nabla}\eta(\xi_{1}{,\xi}_{2}) =
\nabla_{\xi_{1}}\eta( \xi_{2}) - \nabla_{\xi_{2}}\eta(\xi_{1}) -
\eta([\xi_{1}{,\xi}_{2}]).
\]
It is important to note that, in general, $\d_{\nabla}^{2}\neq 0$ so that $\d_{\nabla}$
\emph{is not} a differential. If $\phi:TM\rightarrow A$ is another $A$-valued 1-form,
we define
\[
[\eta,\phi]_{\nabla}(\xi_{1},\xi_{2}) =
T_{\nabla}(\eta(\xi_{1}),\phi(\xi_{2})) +
T_{\nabla}(\phi(\xi_{1}),\eta(\xi _{2}))
\]
where $T_{\nabla}$ is the \emph{torsion} of $\nabla$,
\[T_{\nabla}(\al,\be) = \nabla_{\#\al}\be - \nabla_{\#\be}\al -
[\al,\be]_A, \quad (\al,\be\in\Gamma(A)).\]

We have the following crucial result which follows from a more or less straightforward computation:

\begin{lemma}
An $A$-valued 1-form $\eta$ satisfies the \textbf{generalized Maurer-Cartan equation}:
\begin{equation}
\label{generalizedmaurercartaneq}
d_{\nabla}\eta+\frac{1}{2}[\eta,\eta]_{\nabla} = 0
\end{equation}
if and only if $\eta:TM\to A$ is a Lie algebroid morphism.
\end{lemma}

Now, let $\mathcal{G}$ be a Lie groupoid with Lie algebroid $A$. Right translation by an element
$g\in\G$ is a diffeomorphism between $\mathbf{s}$-fibers: $R_g:\mathbf{s}^{-1}(\mathbf{t}(g))\to \mathbf{s}^{-1}(\mathbf{s}(g))$.
Hence, by a \textbf{right invariant 1-form} on $\mathcal{G}$ we mean a $\mathbf{s}$-foliated $1$-form $\omega$ on $\mathcal{G}$
such that for all $g\in\G$:
\[
\omega(\xi) = \omega(\d_hR_{g}(\xi)), \quad \forall h\in \mathbf{s}^{-1}(\mathbf{t}(g)), \xi\in T_{h}^{\mathbf{s}}\mathcal{G}.
\]
For short, we write this condition as $(R_{g})^{\ast}\omega=\omega$.

A Lie groupoid carries a natural canonical $\mathbf{s}$-foliated right invariant
differential $1$-form with values in its Lie algebroid:

\begin{definition}
The \textbf{Maurer-Cartan form on a Lie groupoid} $\mathcal{G}$ is the $A$-valued
$\mathbf{s}$-foliated right invariant $1$-form defined by
\[
\wmc(\xi) = (\d R_{g^{-1}})_{g}( \xi)
\]
for $\xi\in T_{g}^{\mathbf{s}}\mathcal{G}$.
\end{definition}

Note that the Maurer-Cartan 1-form $\wmc:T^{\mathbf{s}}\G\to A$ covers the target map
$\mathbf{t}:\G\to M$. As one could expect, this 1-form is a solution of the generalized
Maurer-Cartan equation. Moreover, it satisfies a universal property analogous to the case
of Lie groups:

\begin{proposition}
\label{Universal}
Let $\mathcal{G}$ be a Lie groupoid with Lie algebroid $A$
and let $\wmc$ be its right invariant Maurer-Cartan form.
If $\eta:TM\to A$ is a solution of the Maurer-Cartan equation covering
a map $h:M\to X$, then for each $x\in M$ and $g\in\mathcal{G}$ such that
$h(x)=\mathbf{t}(g)$ there exists a unique locally defined diffeomorphism
$\phi:M \to \mathbf{s}^{-1}(h(x))$ satisfying:
\[
\phi(x) = g\text{ and }\phi^{\ast}\wmc = \eta.
\]
\end{proposition}

Let us sketch a proof of this proposition. Since it is a local
result we may assume that $M$ is simply connected. The source 1-connected Lie groupoid
integrating $TM$ is then the pair groupoid $M\times M\tto M$. Since a Maurer-Cartan form
is nothing but a Lie algebroid morphism, we can integrate $\eta:TM\to A$ to a unique Lie
groupoid morphism (see \cite{CrainicFernandes}):
\[
\xymatrix{ M\times M \ar@<.5ex>[d]\ar@<-.5ex>[d]  \ar[r]^-{H} &
\mathcal{G} \ar@<.5ex>[d]\ar@<-.5ex>[d]
\\
M \ar[r]_h & X. }
\]
If we fix a point $x_{0}$ in $M$ we may always write
\[H(x,y) = \phi(x)\phi(y)^{-1}\]
where $\phi:M\rightarrow\mathbf{s}^{-1}(h(x_{0})) \subset\mathcal{G}$ is defined by
$\phi(x):= H(x,x_{0})$. But then $\phi$ satisfies
\[ \phi(x_{0})= \mathbf{1}_{h(x_{0})}, \quad \phi^{\ast}\wmc=\eta. \]
So $\phi$ is the desired local diffeomorphism.
\vskip 10 pt

As a corollary of Proposition \ref{Universal} we obtain:

\begin{corollary}\label{symmetrywmc}
Let $\G\tto X$ be a Lie groupoid with Maurer-Cartan form $\wmc$. If
$\phi : \s^{-1}(x) \rightarrow \s^{-1}(y)$ is a symmetry of $\wmc$
(i.e., $\phi^{\ast}\wmc = \wmc$) then $x$ and $y$ belong to the
same orbit of $\G$ and $\phi$ is locally of the form $\phi = R_g$
for some $g \in \G$.
\end{corollary}

The proof sketched above also shows that if we impose a topological
condition we obtain a global version of the universal
property of Maurer-Cartan forms:

\begin{theorem}
Let $\G\tto X$ be a source 1-connected Lie groupoid with Lie algebroid $A$
and let $\eta\in\Omega^{1}(M,A)$ be an $A$-valued differential 1-form covering $h:M\to X$.
Then, there exists an embedding $\phi:M \to \mathbf{s}^{-1}(h(x_{0}))$
satisfying
\[
\phi(x_0) = 1_{h(x_0)}\text{ and }\phi^{\ast}\wmc = \eta.
\]
if and only if
\begin{enumerate}[(i)]
\item (local obstruction) $\eta$ satisfies the generalized Maurer-Cartan
equation and
\item (global obstruction) the Lie groupoid morphism $H$ integrating $\eta$ is
trivial when restricted to the fundamental group $\pi_1(M,x_{0})$.
\end{enumerate}
\end{theorem}

\subsection{The Classification Theorem}

We can now give a complete solution to the classification problem \ref{PClassification}:

\begin{theorem}
\label{thm:classification}
Let $(n,X,C_{ij}^{k},F_{i}^{a})$ be the initial data of a Cartan's
realization problem. Then:
\begin{enumerate}[(i)]
\item a realization exists if and only if $C_{ij}^k$ and $F_{i}^{a}$ are the structure
functions of a Lie algebroid;
\item any realization is locally equivalent to a neighborhood
of the identity of an $\mathbf{s}$-fiber of a groupoid
$\mathcal{G}$ equipped with the Maurer-Cartan form;
\item two germs of coframes belong to the same global coframe if and only if they
correspond to points on the same orbit of $A$.
\end{enumerate}
\end{theorem}

\begin{proof}
We already know that for Cartan's problem to have a solution the
$C_{ij}^{k}$ and $F_{i}^{a}$ form the structure functions of the
flat Lie algebroid $A=X\times \Rr^n$, with bracket and anchor given
by \eqref{eq:bracket} and \eqref{eq:anchor}. Assume, for simplicity,
that $A$ is integrable and that $\mathcal{G}$ is a Lie groupoid
integrating $A$. Denote by $\wmc$ the Maurer-Cartan form of
$\mathcal{G}$ and by $\{\wmc^{1},\ldots,\wmc^{n}\}$ its components
with respect to the basis $\set{e_1,\ldots,e_n}$. Then it
is clear that, for each $x_{0}\in X$,
$(\mathbf{s}^{-1}(x_{0}),\{\wmc^{i}\},\mathbf{t})$ is a realization of
the Cartan problem with initial data $(n,X,C_{ij}^{k},F_{i}^{a})$. This proves (i).

Next we observe that, if $(M,\{\theta^{i}\},h)$ is another realization of
$(n,X,C_{ij}^{k},F_{i}^{a})$, the $A$-valued $1$-form
\[ \theta=\sum_{i=1}^n \theta^{i}e_i\in\Omega^{1}(M,A)\]
satisfies the generalized Maurer-Cartan equation. In fact, equation \eqref{dh} is equivalent to
$(\theta,h)$ being an $A$-valued differential form and equation \eqref{dtheta} is
equivalent to the Maurer-Cartan equation. Hence, if $p_{0}\in M$ is such that
$h(p_{0})=x_{0}$ then by the universal property of Maurer-Cartan forms (Theorem \ref{Universal})
we can find a neighborhood $V$ of $p_{0}$ in $M$ and a diffeomorphism
\[ \phi:V\to\phi(V)\subset \mathbf{s}^{-1}(x_{0})\]
such that $\phi(p_{0})=\mathbf{1}_{x_{0}}$ and $\phi^{\ast}\wmc=\theta$.
This shows that any realization of Cartan's problem is locally equivalent to
a neighborhood of the identity of an $\mathbf{s}$-fiber of $\mathcal{G}$
equipped with the Maurer-Cartan form, so (ii) follows.

Finally, if two germs of coframes correspond to points on the same orbit of $A$, then it is clear that the restriction of the Maurer-Cartan form to any $\s$-fiber of $\G$ over this orbit gives rise to a coframe containing both germs. We shall now prove the converse. Suppose that $M$ is a connected manifold with a coframe $\theta^i$, and let $x$ and $y$ be points of $M$. The germs of $\theta^i$ at $x$ and at $y$ correspond respectively to the points $h(x)$ and $h(y)$ in $X$. Let $\gamma$ be a curve in $M$ joining $x$ to $y$ and let $\set{U_1, \ldots U_k}$ be an open cover of $\gamma$ with the property that there exists a diffeomorphism $\phi_i$ taking each $U_i$ to some $\s$-fiber and such that $\theta = \phi^{\ast}\wmc$. What we shall show is that if we fix points $x_i \in U_i$ and $x_{i+1} \in U_{i+1}$, then there exists an element $g \in \G$ which maps $h(x_i)$ to $h(x_{i+1})$. For this let $x_{i,i+1} \in U_i \cap U_{i+1}$ and denote by
\[ \begin{array}{lll}
g_i & = & \phi_i(x_i) \\
g_{i+1} & = & \phi_{i+1}(x_{i+1}) \\                                                                                                        g_{i,i+1}& = & \phi_i(x_{i,i+1})  \\                                                                                                       g_{i+1,i}&= & \phi_{i+1}(x_{i,i+1}). \\                                                                                                    \end{array}\]
Then, by construction, we have that
\[ \begin{array}{lll}
\t(g_i) & = & h(x_i) \\
\s(g_{i,i+1})& = & \s(g_i)  \\
\t(g_{i,i+1})& = & \t(g_{i+1,i})  \\
\s(g_{i+1,i})& = & \s(g_{i+1})\\
\t(g_{i+1})& = & h(x_{i+1})\\                                                                                                    \end{array}\]
and thus, $(g_{i+1})(g_{i+1,i})^{-1}(g_{i,i+1})(g_i)^{-1}$ is an element of $\G$ which maps $x_i$ to $x_{i+1}$,
so (iii) also follows.
\end{proof}

Theorem \ref{thm:classification} shows that there is a one
to one correspondence between flat Lie algebroids of rank $n$ and classes of locally defined
coframes on $n$-dimensional manifolds. We call the Lie algebroid associated with
a class of coframes the \textbf{classifying Lie algebroid}. Its integration
furnishes explicit models to the realization problem.

\begin{remark}
Even when $A$ is not integrable, the proof above shows that all we need is a local
groupoid integrating it, and this always exists (see \cite{CrainicFernandes}).
\end{remark}

\subsection{Realizations of Finite Type $G$-Structures}

We now briefly describe the realization problem for $G$-structures
of finite type. We focus on the case where $G^{(1)} = \set{e}$.
There are two reasons for this: (i) this is the case that appears
in most geometric applications, and (ii) more general finite type
$G$-structures can then be handle by induction once this case is
understood.

Suppose that we would like to determine all $G$-structures
belonging to a certain class, where $G$ is a fixed Lie group
satisfying $G^{(1)} = \set{e}$. The considerations of Section
\ref{sec:2ndorder} show that we must solve the following
realization problem:

\begin{problem}[Realization Problem for $G$-Structures with $G^{(1)} = \{e\}$]
Given the data:
\begin{itemize}
  \item An open set $X \subset \Rr^d$,
  \item an integer $n \in \Nn$,
  \item a Lie subalgebra $\gg \subset \mathfrak{gl}_n$ satisfying
  $\gg^{(1)}=0$,
  \item a Lie group $G \subset GL(n)$ with Lie algebra $\gg$, and
  \item maps $c:X \rightarrow \hom (\Rr^n \wedge \Rr^n , \Rr^n)$,
  $b:X \rightarrow \hom ( \Rr^n \wedge \Rr^n , \gg )$,
  $S:X  \rightarrow \hom ( \Rr^n \otimes \gg , \gg)$,
  $\Theta : X \rightarrow \hom ( \Rr^n , \Rr^d )$, and
  $ \Phi : X \rightarrow \hom ( \gg , \Rr^d )$
\end{itemize}
one asks for the existence of
\begin{enumerate}[(1)]
  \item a manifold $M^n$,
  \item a $G$-structure $\BG(M)$ on $M$ with tautological form
  $\omega \in \Omega^1(\BG , \Rr^n)$
  \item a maximal rank one form $\phi \in \Omega^1(\BG, \gg)$
  transversal to $\omega$, and
  \item a map $h: \BG \rightarrow X$
\end{enumerate}
such that:
\begin{align}
\d\omega &= c(h)\circ\omega\wedge\omega - \phi\wedge\omega \\
\d\phi   &= b(h)\circ\omega\wedge\omega + S(h)\circ\omega\wedge\phi - \phi\wedge\phi \\
\d h     &= \Theta(h)\circ\omega + \Phi(h)\circ\phi
\end{align}
\end{problem}

Just like in the case of $G=\{e\}$, one checks that a necessary condition
for solving this problem is that the structure functions $c,b,S,\Theta,\Phi$
determine a flat Lie algebroid $A \rightarrow X$ with fiber $\Rr^n \oplus
\gg$. The Lie bracket on constant sections, is given by
\[
[(u,\al),(v,\be)]_A(x)=(w,\ga),
\]
where,
\begin{align*}
w   &:=c(x)(u \wedge v)-\al\cdot v + \be\cdot u, \\
\ga &:=b(x)(u \wedge v)+S(x)(u\otimes\be-v\otimes\al) - [\al ,\be]_{\gg}.
\end{align*}
The anchor map $\#:A\to TX$ is determined by:
\[
\#(u,\al):=\Theta(x)u+\Phi(x)\al.
\]

Observe that the natural inclusion $\gg\hookrightarrow \Gamma(\A)$, $\al\mapsto(0,\al)$,
is a Lie algebra homomorphism, since we have:
\[
[(0,\al),(0,\be)]_A = (0,[\al ,\be]_{\gg}).
\]
Therefore, we obtain an inner action of $\gg$ on $\A$ by setting
\[
\rho:\gg\to\Der(A),\quad \rho(\al)(\sigma):=[(0,\al),\sigma]_{\A} \quad (\sigma\in\Gamma(A)).
\]
Conversely, if the initial data of the realization problem
determines a flat Lie algebroid $A=X\times(\Rr^n\oplus\gg) \rightarrow X$
such that the natural inclusion $\gg\hookrightarrow \Gamma(\A)$ defines
an infinitesimal inner $\gg$-action, then for each $h_0\in X$ we
can find a $G$-structure whose structure functions assume the
values $c(h_0),b(h_0),S(h_0),\Theta(h_0)$ and $\Phi(h_0)$.

Note that now, in order to exhibit explicit universal models, we must
integrate not only the Lie algebroid but also the infinitesimal
inner $\gg$-action. Therefore, assume that $A$ is integrable and that
we can integrate the $\gg$-action to a proper and free action of $G$ on
some groupoid $\G$ integrating $A$. Each $\s$-fibers of $\G$ will be
the total space of a $G$-structure determined by its tautological form; the $\Rr^n$ component of the
Maurer-Cartan form. The isotropy group of $\G$ at a point will be
the symmetry group of the corresponding $G$-structure. Moreover,
for any frame $p$ in another realization $\BG(M)$, there will
exist a neighborhood which is isomorphic to a neighborhood of the
identity in an $\s$-fiber of $\G$. Details on these constructions will
appear elsewhere.

\section{Examples}

In order to illustrate our method, we will present two concrete examples.
The first one is a toy example: metrics of constant curvature in
$\Rr^2$. The second example is a more serious application to
Bochner-K\"{a}hler metrics and is based on the work of Bryant \cite{Bryant}.

\subsection{Metrics of Constant Curvature in $\Rr^2$}
Suppose we would like to classify all Riemann metrics of constant
curvature in a neighborhood of the origin in $\Rr^2$. The
$G$-structure to be considered is $\B_{O_2}(\Rr^2)$, the orthogonal frame bundle of
$\Rr^{2}$. In this case, the first order structure function vanishes identically and
$O_{2}^{(1)}= \{e\}$. The structure equation has the following
simple form:
\[
\left\{
\begin{array}{l}
\d \omega = -\phi\wedge\omega\\
\d \phi   = k\omega\wedge\omega,
\end{array}
\right.
\]
where $\phi$ is the connection 1-form on $\B_{O_2}(\Rr^2)$ corresponding to the
Levi-Civita connection and $k$ is the Gaussian curvature of
$\phi$. If the curvature $k$ is constant, then $\d k=0$. In terms of
the canonical base of $\Rr^{2}$ the structure equation becomes
\begin{align*}
\d \omega^{1} &= -\phi\wedge\omega^{2},\\
\d \omega^{2} &= \phi\wedge\omega^{1},\\
\d \phi       &= k\omega^{1}\wedge\omega^{2},\\
\d k          &= 0.
\end{align*}
The Gaussian curvature is the only invariant function. It follows
that the classifying Lie algebroid for this class is the flat bundle
$A=\Rr\times\Rr^3\rightarrow \Rr$, with a basis of sections
$\{e_{1},e_{2},e_{3}\}$ and structure given by:
\begin{align*}
\left[ e_{1} , e_{2} \right] (k) &= ke_{3}, \\
\left[ e_{1} , e_{3} \right] (k) &= -e_{2}, \\
\left[ e_{2} , e_{3} \right](k)  &= e_{1}, \\
\# e_{i}  &= 0.
\end{align*}
Note that this Lie algebroid is just a bundle of Lie algebras with fibers isomorphic to:
\begin{center}
\begin{tabular}{||c|c|c||}\hline
$\mathfrak{sl}_{2}$   & \footnotesize{if $k < 0$} & \footnotesize{Hyperbolic Geometry}  \\
\hline $\mathfrak{se}_{2}$ & \footnotesize{if $k=0$}   & \footnotesize{Euclidean Geometry}        \\
\hline $\mathfrak{so}_{3}$ & \footnotesize{if $k>0$} & \footnotesize{Spherical Geometry}  \\
\hline
\end{tabular}
\end{center}
\vskip 10 pt
The inner action of $\mathfrak{o}_{2}$ on $A$ is the fiberwise
action obtained from the adjoint representation. It follows that
to each value of $k\in\Rr$ there corresponds the germ at
$0$ of a constant curvature metric on $\Rr^{2}$. Moreover,
two constant curvature metrics on $\Rr^{2}$ are locally equivalent
if and only if they have the same curvature.

\subsection{Bochner-K\"{a}hler Structures}
Let $(M,\sigma,\Omega)$ be a K\"{a}hler manifold. Its curvature tensor
can be decomposed into three irreducible components: the scalar curvature,
the traceless Ricci curvature and the Bochner curvature. We say that
$(M,\sigma,\Omega)$ is Bochner-K\"{a}hler if its Bochner curvature vanishes
identically.

Bryant in \cite{Bryant} performs a differential analysis which shows that the
local classification of Bochner-K\"{a}hler metrics can be reduced to the
following structure equations
\begin{align*}
\d\omega &= -\phi\wedge\omega, \\
\d\phi   &= -\phi\wedge\phi + S\omega^{\ast}\wedge\omega-S\omega\wedge\omega^{\ast} - \omega\wedge\omega^{\ast}S +
                (\omega^{\ast}\wedge S\omega)I_{n}, \\
\d S     &= -\phi S+S\phi+T\omega^{\ast}+\omega T^{\ast} +
                \frac{1}{2}(T^{\ast}\omega+\omega^{\ast}T)I_{n}, \\
\d T       &= -\phi T + (UI_{n}+S^{2}) \omega, \\
\d U     &= T^{\ast}S\omega+\omega^{\ast}ST.
\end{align*}
where $\omega$ is the $\mathbb{C}^{n}$-valued tautological form on
the unitary coframe bundle of $M$, $\phi$ is the
$\mathfrak{u}(n)$-valued connection form associated to the
Levi-Civita connection on $M$, $S$,$T$, and $U$ are functions with
values, respectively, in $i\mathfrak{u}(n)$, $\mathbb{C}^{n}$ and $\Rr$
and $I_{n}\in GL_{n}(\mathbb{C})$ is the identity.

The corresponding classifying Lie algebroid can be described as follows: as a
vector bundle, it is the trivial bundle over $X = i\mathfrak{u}(n)
\oplus \mathbb{C}^{n} \oplus \Rr$ with fiber type
$\mathbb{C}^{n} \oplus \mathfrak{u}(n)$. The Lie bracket is defined
on constant sections by
\begin{multline*}
[(z_{1},\alpha_{1}),(z_{2},\alpha_{2})]_A(s,t,u) :=
(\alpha_{2}z_{1}-\alpha_{1}z_{2}, -
[\alpha_{1},\alpha_{2}]_{\mathfrak{u}(n)} + (z_{1}^{\ast}z_{2} -
z_{2}^{\ast}z_{1})s + \\
- s (z_{1}z_{2}^{\ast} - z_{2}z_{1}^{\ast}) -
(z_{1}(z_{2}^{\ast}s) - z_{2}(z_{1}^{\ast}s)) +
(z_{1}^{\ast}(sz_{2}) - z_{2}^{\ast}(sz_{1}))I_{n}),
\end{multline*}
while the anchor is given by:
\begin{multline*}
\# (z,\alpha) (s,t,u) := (-\alpha s + s\alpha + tz^{\ast} + zt^{\ast} +
\frac{1}{2}(t^{\ast} z + z^{\ast}t) I_{n} , \\
-\alpha t + (uI_{n}+s^{2}) z , t^{\ast}sz + z^{\ast}st).
\end{multline*}
Finally, the inner action of $\mathfrak{u}(n)$ on $A=X\times(\mathbb{C}^{n} \oplus \mathfrak{u}(n))$ is composed of the defining
action of $\mathfrak{u}(n)$  on $\mathbb{C}^{n}$ and the adjoint action on $\mathfrak{u}(n)$. This action covers an infinitesimal
$\mathfrak{u}(n)$-action on $X=i\mathfrak{u}(n)\oplus \mathbb{C}^{n} \oplus \Rr$, which on the first factor corresponds to the
adjoint action, on the second factor is just the defining action of $\mathfrak{u}(n)$ and is trivial on the last factor.

In order to understand the local classification of Bochner-K\"{a}hler metrics and find local models for such metrics
one must integrate this classifying Lie algebroid. The orbit space will correspond to the different classes of Bochner-K\"{a}hler metrics,
while the different isotropy types will correspond to their symmetry groups. We leave this analysis
for a future work.

\bibliographystyle{amsplain}
\def\lllll{}

\end{document}